\documentclass{article}
\usepackage{yfonts}
\usepackage{amsmath,amssymb,authblk,amsthm}
\usepackage{graphicx}
\usepackage[all]{xy}
\def\<{\left <}
\def\rra{\rightrightarrows}
\def\>{\right >}
\def\({\left (}
\def\){\right )}
\def\bra {\left [}
\def\ket{\right ]}
\def\cbra{\left \{}
\def\cket{\right \}}

\DeclareSymbolFont{AMSb}{U}{msb}{m}{n}
\DeclareMathSymbol{\N}{\mathbin}{AMSb}{"4E}
\DeclareMathSymbol{\Z}{\mathbin}{AMSb}{"5A}
\DeclareMathSymbol{\R}{\mathbin}{AMSb}{"52}
\DeclareMathSymbol{\Q}{\mathbin}{AMSb}{"51}
\DeclareMathSymbol{\I}{\mathbin}{AMSb}{"49}
\DeclareMathSymbol{\C}{\mathbin}{AMSb}{"43}

\newtheorem{thm}{Theorem}[section]
\newtheorem{defn}[thm]{Definition}
\newtheorem{lem}[thm]{Lemma}
\newtheorem{prop}[thm]{Proposition}

\begin{document}

\date{}
\title{ Non-deterministic Semantics for Dynamic Topological Logic
      }
\author{David Fern\'{a}ndez-Duque}
\affil{Department of Mathematics, Stanford University, Stanford, CA
94305, USA
}

\maketitle

\begin{abstract}
Dynamic Topological Logic ($\mathcal{DTL}$) is a combination of $\mathcal{S}${\em 4}, under its topological interpretation, and the temporal logic $\mathcal{LTL}$ interpreted over the natural numbers. $\mathcal{DTL}$ is used to reason about properties of dynamical systems based on topological spaces. Semantics are given by dynamic topological models, which are tuples $\left <X,\mathcal{T},f,V\right >$, where $\left <X,\mathcal{T}\right >$ is a topological space, $f$ a function on $X$ and $V$ a truth valuation assigning subsets of $X$ to propositional variables.

Our main result is that the set of valid formulas of $\mathcal{DTL}$ over spaces with continuous functions is recursively enumerable. We do this by defining alternative semantics for $\mathcal{DTL}$. Under the standard semantics, $\mathcal{DTL}$ is not complete for Kripke frames. However, we introduce the notion of a non-deterministic quasimodel, where the function $f$ is replaced by a binary relation $g$ assigning to each world multiple temporal successors. We place restrictions on the successors so that the logic remains unchanged; under these alternative semantics, $\mathcal{DTL}$ becomes Kripke-complete. We then apply model-search techniques to enumerate the set of all valid formulas.
\end{abstract}

\section{Introduction}
{\em Dynamic Topological Logic} ($\mathcal{DTL}$) is a propositional tri-modal system introduced in \cite{kmints,arte} for reasoning about topological dynamics; that is, about the action of a continuous function $f$ on a topological space $\<X,\mathcal T\>$. The interpretation of formulas of $\mathcal{DTL}$ involves not only points $x\in X$, but also their orbit
\[\cbra x,f(x),f^2(x),...\cket\]
and their neighborhoods in $\mathcal T$.  

The language uses propositional variables, Classical Boolean connectives and three modalities: $\Box$ from $\mathcal{S}${\em 4}, interpreted as the topological interior, and the temporal operators $\bigcirc$ and $\ast$ of Linear Temporal Logic (\cite{temporal}), which are interpreted as `next' and `henceforth', respectively.

Every class $\mathcal{C}$ of dynamic topological systems induces a logic consisting of those formulas of $\mathcal{DTL}$ that are valid on all systems in $\mathcal{C}$. Many logics arising in this form have been studied; the following are some of the main results which are known.

\begin{enumerate}
\item The fragment $\mathcal{DTL}^\bigcirc$:

$\mathcal{DTL}^\bigcirc$ is the fragment of $\mathcal{DTL}$ which uses only $\Box$ and $\bigcirc$. This fragment is finitely axiomatizable and has the finite model property, both when $f$ is taken to be a continuous function (\cite{arte}) and when it is a homeomorphism (\cite{kmints}). Interestingly enough, the logic over arbitrary spaces does not coincide with the logic over $\R$ (\cite{slav2}) but it does coincide with the logic over $\Q$ (\cite{rationals}).

\item The fragment $\mathcal{DTL}_0$:

In this fragment all three modalities are used but temporal modalities may not appear in the scope of $\Box$. $\mathcal{DTL}_0$ is finitely axiomatizable and decidable (\cite{kmints}).

\item The fragment $\mathcal{DTL}_1$:

$\mathcal{DTL}_1$ is the fragment of $\mathcal{DTL}$ where all three modalities are used but $\ast$ may not appear in the scope of $\Box$. This fragment is powerful enough to encode some undecidable problems and hence is undecidable. However, it is recursively enumerable (\cite{konev}), and the logic coincides over arbitrary spaces, locally finite spaces and $\R^2$ (\cite{konev,slav,me}). This implies that all these logics are also equal on the smaller fragments $\mathcal{DTL}^\bigcirc$ and $\mathcal{DTL}_0$.

\item Spaces with homeomorphisms:

The set of valid formulas of $\mathcal{DTL}_1$ over spaces with homeomorphisms is not recursively enumerable (\cite{wolter}), and hence the same is true for the set of valid formulas of the full language. Furthermore, the logics over arbitrary spaces, Aleksandroff spaces and $\R^n$ for $n\geq 0$ are all distinct (\cite{wolter,slav2}).

\item Full $\mathcal{DTL}$ with arbitrary continuous functions:

It is known that for full $\mathcal{DTL},$ the logics over arbitrary spaces, Aleksandroff spaces and $\R^n$ are all distinct (\cite{kmints,me}). Over almost disjoint spaces (where all open sets are closed) the logic is decidable, even when $f$ is taken to be a homeomorphism (\cite{s5}).

$\mathcal{DTL}$ is over arbitrary spaces is undecidable; this follows from the fact that $\mathcal{DTL}_1$ is already undecidable. However, it is recursively enumerable; this is the main result we will present here.
\end{enumerate}

The layout of this paper is as follows. In \S\ref{secndq} we will define non-deterministic quasimodels for $\mathcal{DTL}$, where the function of a dynamic topological system is replaced by a binary relation $g$, so that each point may have several immediate temporal successors. However, we will place restrictions on $g$ so that the logic remains sound, and in \S\ref{secgen} show that a dynamic topological model can always be reconstructed from a non-deterministic quasimodel.

The central result is that $\mathcal{DTL}$ is complete for the class of locally finite Kripke frames under the new semantics. Our strategy for proving this will be to generate truth-preserving binary relations between Kripke frames and dynamic topological models. The relations we will use are called {\em $\omega$-simulations} and are developed in \S\ref{secsim}.

We then apply techniques very similar to those in \cite{konev} to show that $\mathcal{DTL}$ is recursively enumerable. There, Kruskal's Tree Theorem is used to prove that a certain model-search algorithm always reports failure in finite time when a non-satisfiable formula of $\mathcal{DTL}_1$ is given as imput. In \S\ref{re} we use non-deterministic semantics to develop a variation of this which can be applied to arbitrary formulas of the language.
\section{Dynamic Topological Logic}\label{basic}

The language of $\mathcal{DTL}$ is built from propositional variables in a countably infinite set ${\bf Var}$ using the Boolean connectives $\wedge$ and $\neg$ (all other connectives are to be defined in terms of these) and the three unary modal operators $\Box$ (`interior'), $\bigcirc$ (`next') and $\ast$ (`henceforth'). We write $\Diamond$ as a shorthand for $\neg\Box\neg$. Formulas of this language are interpreted on dynamical systems over topological spaces, or {\em dynamic topological systems}.

\begin{defn}
A {\em dynamic topological system} is a triple $\mathfrak S=\<X,\mathcal{T},f\>,$
where $\<X,\mathcal{T}\>$ is a topological space and $f:X\to X$ is a continuous function.

A {\em valuation on $\mathfrak S$} is a relation $V\subseteq{\bf Var}\times X.$ A dynamic topological system equipped with a valuation is a {\em dynamic topological model}.
\end{defn}

The valuation $V$ is extended inductively to arbitrary formulas as follows:
\[
\begin{array}{lcl}
V(\alpha\wedge\beta)&=&V(\alpha)\cap V(\beta)\\
V(\neg\alpha)&=&X\setminus V(\alpha)\\
V(\Box\alpha)&=&V(\alpha)^\circ\\
V(\bigcirc\alpha)&=&f^{-1}V(\alpha)\\
V(\ast\alpha)&=&\displaystyle\bigcap_{n\geq 0}f^{-n}V(\alpha).
\end{array}
\]

$\mathcal{DTL}$ distinguishes arbitrary spaces from finite spaces and even from locally finite spaces (those where every point has a neighborhood with finitely many points).

More generally, $\mathcal{DTL}$ distinguishes arbitrary topological spaces from Aleksandroff spaces (\cite{kmints}); that is, spaces where arbitrary intersections of open sets are open (\cite{aleksandroff}). All locally finite spaces are Aleksandroff spaces.

Nevertheless, we will show how locally finite spaces can be used to represent a larger class dynamic topological systems; to do this, we will define alternative semantics for $\mathcal{DTL}$.

\section{Non-deterministic quasimodels}\label{secndq}
We will denote the set of subformulas of $\varphi$ by ${\mathrm{sub}}(\varphi)$, and define
$${\mathrm{sub}}_{\pm}(\varphi)={\mathrm{sub}}(\varphi)\cup\neg {\mathrm{sub}}(\varphi).$$
If we identify $\psi$ with $\neg\neg\psi$, one can think of ${\mathrm{sub}}_{\pm}(\varphi)$ as being closed under negation.

A set of formulas ${\bf t}\subseteq {\mathrm{sub}}_\pm(\varphi)$ is a {\em $\varphi$-type} if, for all $\psi\in {\mathrm{sub}}_\pm(\varphi),$
\[\psi\not\in{\bf t}\Leftrightarrow \neg\psi\in {\bf t}\]
and for all $\psi_1\wedge\psi_2\in {\mathrm{sub}}_\pm(\varphi),$
$$\psi_1\wedge\psi_2\in {\bf t}\Leftrightarrow \psi_1\in {\bf t}\text{ and }\psi_2\in{\bf t}.$$

The set of $\varphi$-types will be denoted by ${\bf type}(\varphi)$.

\begin{defn}[typed Kripke frame]\label{frame}
Let $\varphi$ be a formula in the language of $\mathcal{DTL}$.
A {\em $\varphi$-typed Kripke frame} is a triple
$\mathfrak{F}=\<W,R,t\>,$
where $W$ is a set, $R$ a transitive, reflexive relation on $W$ and $t$  a function assigning a $\varphi$-type $t(w)$ to each $w\in W$ such that
\[\Box\psi\in t(w)\Leftrightarrow \forall v\(Rwv\Rightarrow\psi\in t(v)\).\]
\end{defn}

It is easy to see that this is equivalent to the dual condition that
\[\Diamond\psi\in t(w)\Leftrightarrow \exists v\(Rwv\text{ and }\psi\in t(v)\).\]
Kripke frames give sound and complete semantics for for $\mathcal{S}${\em 4} (\cite{black}), but here we are disregarding the temporal modalities by giving valuations of these formulas a priori rather than by their usual meaning. One would then be tempted to equip the Kripke frame with a transition function in order to interpret temporal operators directly. However, this would give us a class of models for which $\mathcal{DTL}$ is incomplete; instead, we will allow each world to have multiple temporal successors via a `sensible' relation $g$, as we will define below. For our purposes, a {\em continuous relation} on a topological space is a relation under which the preimage of any open set is open.

\begin{defn}[sensible relation]\label{compatible}
Let $\varphi$ be a formula of $\mathcal{DTL}$ and $\<W,R,t\>$ a $\varphi$-typed Kripke frame.

Suppose that ${\bf t},{\bf s}\in{\bf type}(\varphi)$. The ordered pair $({\bf t},{\bf s})$ is {\em sensible} if
\begin{enumerate}
\item for all $\bigcirc\psi\in {\mathrm{sub}}(\varphi)$,
$\bigcirc\psi\in {\bf t}\Leftrightarrow \psi\in {\bf s}$ and
\item for all $\ast\psi\in {\mathrm{sub}}(\varphi)$,
$\ast\psi\in {\bf t}\Leftrightarrow\(\psi\in {\bf t}\text{ and }\ast\psi\in {\bf s}\).$
\end{enumerate}

Likewise, a pair $(w,v)$ of worlds in $W$ is sensible if $(t(w),t(v))$ is sensible.

A continuous relation
$g\subseteq W\times W$
such that $g(w)\not=\varnothing$ for all $w\in W$ is {\em sensible} if every pair in $g$ is sensible.

Further, $g$ is $\omega$-sensible if for all $\ast\psi\in {\mathrm{sub}}(\varphi)$,
\[\neg\ast\psi\in t(w)\Leftrightarrow\exists v\in W\text{ and }N\geq 0\text{ such that }\neg\psi\in t(v)\text{ and }g^Nwv.\]
\end{defn}

It is a good idea to examine what continuity means in a Kripke frame. Suppose $\mathfrak{F}$ is as in Definition \ref{frame} and $g$ is continuous. Pick $w,v\in W$ so that $gwv$. Since the set
$R(v)=\cbra u:Rvu\cket$
is open, we know that
$R(w)\subseteq g^{-1}R(v).$
In other words, if $Rww'$, there exists $v'$ such that $Rvv'$ and $gw'v',$ so that the following square can always be completed:
$$
\xymatrix
{
w'\ar@{-->}[r]^g&v'\\
w\ar[u]^R\ar[r]^g&v\ar@{-->}[u]^R.
}
$$

We are now ready to define our non-deterministic semantics for $\mathcal{DTL}$.

\begin{defn}[non-deterministic quasimodel]\label{ndqm}
A {\em $\varphi$-typed non-deterministic quasimodel} is a tuple
$\mathfrak{ D}=\<W,R,t,g\>,$
where $\<W,R,t\>$ is a $\varphi$-typed Kripke frame and $g$ is an $\omega$-sensible relation on $W$.

$\mathfrak D$ {\em satisfies} $\varphi$ if there exists $w_\ast\in W$ such that $\varphi\in t(w_\ast)$.
\end{defn}

Non-deterministic quasimodels are similar to dynamic Kripke frames (\cite{arte,me,konev}) except for the fact that $g$ is now a relation instead of a function. However, note that we do not allow any subformulas of $\varphi$ to be left undecided by $g$; that is, if $\bigcirc\psi\in {\mathrm{sub}}_{\pm}(\varphi)$ and $w\in W$, then either $\psi\in t(v)$ for all temporal successors $v$ of $w$ or $\neg\psi\in t(v)$ for all such $v$. This is necessary in order to preserve soundness.

\begin{figure}[htp]
\begin{center}
\scalebox
{0.9}
{
$$
\xymatrix{
v \ar@(u,l)[]_g\ar[r]^g & w\ar@(r,u)[]_g\\
u \ar@(d,l)[]^g\ar@{~>}[u]^R & \\
}
$$
}
\end{center}
\caption
{
The formula $\varphi=\ast\Box p\to\Box\ast p$ is valid on all finite topological models; in fact, it is valid on all topological models based on a locally finite topological space (\cite{kmints}). However, $\varphi$ can be refuted in a non-deterministic quasimodel with only three worlds. Take $W=\cbra u,v,w\cket$ and let $R,g$ be as shown in the diagram above (closing $R$ under reflexivity). Assign types to $u,v$ and $w$ in such a way that $p,\ast\Box p,\neg\Box\ast p\in t(u)$, $p,\neg\ast p\in t(v)$ and $\neg p\in t(w)$. Then, $\varphi\not\in t(u)$, hence we have a non-deterministic quasimodel satisfying $\neg\varphi$. This, we will see below, shows that $\varphi$ is not a theorem of $\mathcal{DTL}$.
}
\label{example}
\end{figure}
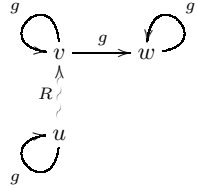

\section{Generating dynamic topological models from non-de\-ter\-mi\-nis\-tic quasimodels}\label{secgen}

A dynamic topological model can be constructed from any non-deterministic quasimodel. Evidently, a non-deterministic quasimodel is not always a dynamic topological model, since dynamic topological models require a transition function rather than a relation. Instead, we will build a topological space whose points are infinite sequences of worlds.

\subsection{Realizing sequences}
A {\em path} in $\mathfrak D$ is any finite or infinite sequence $\<w_n\>$ such that $gw_nw_{n+1}$.

An infinite path
${\vec w}=\< w_n\>_{n\geq 0}$
is {\em realizing} if for all $n\geq 0$ and $\neg\ast\psi\in t(w_n)$ there exists $K\geq n$ such that $\neg\psi\in t(w_K)$.

Denote the set of realizing paths by $W^g$. We will construct dynamic topological models from non-deterministic quasimodels by topologizing this set. All other paths will be thrown away, since in such paths, $\ast$ could not be interpreted according to its intended meaning.

The main transformation we will consider on $W^g$ will be the `shift' operator, defined by
$\sigma\(\<w_n\>_{n\geq 0}\)=\<w_{n+1}\>_{n\geq 0}.$
This simply removes the first element in the sequence.

For our construction to work we must guarantee that there are `enough' realizing paths, in the sense of the following definition.
\begin{defn}[extensive subset]
Let $\varphi$ be a formula of $\mathcal{DTL}$,
\[\mathfrak{D}=\<W,R,t,g\>\]
be a $\varphi$-typed non-deterministic quasimodel and $Y\subseteq W^g$.

Then, $Y$ is {\em extensive} if
\begin{enumerate}
\item $Y$ is closed under $\sigma$;
\item any finite path
$\<w_0,w_1,..,w_N\>$
in $\mathfrak D$ can be extended to an infinite path
$
{\vec w}=\< w_n\>_{n\geq 0}\in Y.
$
\end{enumerate}
\end{defn}
\begin{lem}\label{extension}
If
$\mathfrak D=\<W,R,t,g\>$
is a non-deterministic quasimodel, $W^g$ is extensive.
\end{lem}

\proof
It is obvious that $W^g$ is closed under $\sigma$.

Let
$\<w_0,...,w_N\>$
be a finite path and
$\psi_0,...,\psi_I$
be all formulas such that $\neg\ast\psi_i\in t(w_N)$. Because $g$ is $\omega$-sensible, we know that there exist $K_I$ and $v_I\in g^{K_I}(w_N)$ such that $\neg\psi_I\in t\(v_I\)$. We can then define $w_{N+1},...,w_{N+K_I}=v_I$
in such a way that $gw_nw_{n+1}$ for all $n<N+K_I.$
Now consider $\psi_{I-1}$. If
$\neg\psi_{I-1}\in t(w_n)$
for some $n\leq K_I,$ there is nothing to do and we can set $K_{I-1}=0$.

Otherwise,
$\neg\ast\psi_{I-1}\in t(w_{N+K_I})$
and we can pick $K_{I-1}$ and $v_{I-1}$ such that
$g^{K_{I-1}}v_{I-1}v_{I}.$
Then, define
$w_{N+K_I+1},...,w_{N+K_I+K_{I+1}}=v_{I-1}$
as before.

Continuing inductively, we can define $\cbra w_n\cket_{n\leq N+K},$ where
$K=\sum_{i\leq I}K_i$
and for all $I\leq i,$ $\neg\psi_i\in t(w_{N+k})$ for some $k\leq K$. We can then repeat the process starting with $\cbra w_n\cket_{n\leq N+K},$ and continue countably many times to get a path $\cbra w_n\cket_{n\geq 0}.$ It is then easy to see that this path is realizing.
\footnote{
We must ensure that $K>0$ at each step so that the sequence increases in length and the end result is an infinite path, but this can always be done since $g(w)\not=\varnothing$ for all $w$. 
}
\endproof
\begin{lem}\label{pathcont}
Let
$\mathfrak{D}=\<W,R,t,g\>$
be a $\varphi$-typed non-deterministic quasimodel,
$\<w_n\>_{n\leq N}$
a finite path and $v_0$ be such that $Rw_0v_0.$

Then, there exists a path
$\<v_n\>_{n\leq N}$
such that, for $n\leq N$, $Rw_nv_n$.
\end{lem}
\proof
This follows from continuity of $g$ by an easy induction on $N$.
\endproof
\subsection{Limit models}

If $\varphi$ is a formula of $\mathcal{L}$ and
$\mathfrak{D}=\<W,R,g,t\>$
is a non-deterministic quasimodel, the relation $R$ induces a topology on $W$, as we have seen before, by letting open sets be those which are upward closed under $R$. Likewise, $R$ induces a very different topology on $W^g$, in a rather natural way:

\begin{lem}\label{basis}
For each ${\vec w}\in W^g$ and $N\geq 0$
define
\[R_N\({\vec w}\)=\cbra\cbra v_n\cket_{n\geq 0}\in W^g:\forall n\leq N, Rw_nv_n\cket.\]

Then, the set
$\mathcal B^R=\cbra R_N\({\vec w}\):{\vec w}\in W^g,N\geq 0\cket$
forms a topological basis on $W^g$.
\end{lem}
\proof Recall that a collection $\mathcal{B}$ of subsets of $X$ is a basis if
\begin{enumerate}
\item $\displaystyle\bigcup_{B\in\mathcal{B}}B=X$;
\item whenever $B_1,B_2\in\mathcal{B}$ and $x\in B_1\cap B_2$, there exists $B_3\subseteq B_1\cap B_2$ such that $x\in B_3.$
\end{enumerate}

To check the first property, note that it is obvious that, given any path ${\vec w}\in W^g,$ there is a basic set containing it (namely, $R_0\(\vec w\)$). Hence
$W^g=\bigcup_{{\vec w}\in W^g}R_0\(\vec w\).$

As for the second, assuming that ${\vec w}\in R_{N_0}\({\vec w}_0\)\cap R_{N_1}\({\vec w}_1\),$ one can see that
$R_{\max\(N_0,N_1\)}\({\vec w}\)\subseteq R_{N_0}\({\vec w}_0\)\cap R_{N_1}\({\vec w}_1\)$
using the transitivity of $R$.\endproof

\begin{defn}
The topology $\mathcal{T}^R$ on $W^g$ is the topology generated by the basis $\mathcal B^R$.
\end{defn}

Now that we have equipped $W^g$ with a topology, we need a continuous transition function on it to have a dynamic topological system.

\begin{lem}
The `shift' operator
$\sigma:W^g\to W^g$
is continuous under the topology $\mathcal{T}^R$.
\end{lem}
\proof Let
${\vec w}=\cbra w_n\cket_{n\geq 0}$
be a realizing path and $R_N\(\sigma\({\vec w}\)\)$ be a neighborhood of $\sigma\({\vec w}\)$. Then, if
${\vec v}\in R_{N+1}\({\vec w}\),$
$Rw_nv_n$ for all $n\leq N+1$, so $Rw_{n+1}v_{n+1}$ for all $n\leq N$ and
$\sigma\(\vec v\)\in R_N\(\sigma\(\vec w\)\).$
Hence
$\sigma\(R_{N+1}\({\vec w}\)\)\subseteq R_N\(\sigma\(\vec w\)\),$
and $\sigma$ is continuous.\endproof

Finally, we will use $t$ to define a truth valuation: if $p$ is a propositional variable, set
$V^t(p)=\cbra{\vec w}\in W^g:p\in t\(w_0\)\cket.$

We are now ready to assign a dynamic topological model to every non-deterministic quasimodel:

\begin{defn}[limit model]
Given a non-deterministic quasimodel
$\mathfrak{D}=\<W,R,g,t\>,$
define
$\lim \mathfrak{D}=\<W^g,\mathcal{T}^R,\sigma,V^t\>$
to be the {\em limit model} of $\mathfrak{D}$.
\end{defn}

Of course this model is only useful if $V^t$ corresponds with $t$ on all subformulas of $\varphi$, not just propositional variables. Fortunately, this turns out to be the case.

\begin{lem}\label{sound}
Let $Y\subseteq W^g$ be extensive,
${\vec w}=\cbra w_n\cket_{n\geq 0}\in Y$
and $\psi\in {\mathrm{sub}}_{\pm}(\varphi)$.
Then,
\[\<\lim \mathfrak D\upharpoonright Y,{\vec w}\>\models \psi\text{ if and only if } \psi\in t(w_0).\]
\end{lem}
\proof The proof goes by standard induction of formulas. The induction steps for Boolean operators are trivial; here we will only treat the cases for the modal operators.

\paragraph{Case 1: $\psi=\Box\alpha.$}
If $\psi\in t\(w_0\)$, take the neighborhood $R_0\(w_0\)$ of ${\vec w}.$ We can then see that
$$
\begin{array}{rcl}
\Box\alpha\in t(w_0)&\Rightarrow&\forall v \(Rw_0v\Rightarrow \alpha\in t(v)\)\\
&\Rightarrow&\forall {\vec v}\in R_0({\vec w}),\alpha\in t(v_0)\\
\text{\small IH}&\Rightarrow&\forall {\vec v}\in R_0({\vec w}),\<\lim \mathfrak D\upharpoonright Y,{\vec v}\>\models\alpha\\
&\Rightarrow&\<\lim \mathfrak D\upharpoonright Y,{\vec w}\>\models\Box\alpha.\\
\end{array}
$$

On the other hand, if $\psi\not\in t(w_0)$, any neighborhood ${U}_{{\vec w}}$ of ${\vec w}$ contains a subneighborhood $R_{N}({\vec w})$ for some $N\geq 0$ (because these sets generate the topology). Then, by Lemma \ref{pathcont}, there exists a path
$\<v_0,...,v_N\>\subseteq W$
such that $Rw_nv_n$, $gv_nv_{n+1}$ and $\neg\alpha\in t(v_0)$. Because $Y$ is extensive,
$\cbra v_n\cket_{0\leq n\leq N}$
can be extended to a realizing path ${\vec v}\in Y$. Then ${\vec v}\in U_{{\vec w}}$, and by induction hypothesis we have that
$\<\lim \mathfrak D\upharpoonright Y,{\vec v}\>\models\neg\alpha.$

Since ${U}_{{\vec w}}$ was arbitrary, we conclude that
$\<\lim \mathfrak D\upharpoonright Y,{\vec w}\>\models \neg\Box\alpha.$

\paragraph{Case 2: $\psi=\bigcirc\alpha.$} This case follows from the fact that $(w_0,w_1)$ is sensible.

\paragraph{Case 3: $\psi=\ast\alpha.$} 
Because ${{\vec w}}$ is a realizing path, we have that if $\neg\ast\alpha\in t\(w_0\),$ $\neg\alpha\in t(w_N)$ for some $N\geq 0$. We can use the induction hypothesis to conclude that
$\<\lim \mathfrak D\upharpoonright Y,{\sigma}^N\(\vec w\)\>\not\models \alpha$
and so
\[\<\lim \mathfrak D\upharpoonright Y,{{\vec w}}\>\not\models\ast\alpha.\]

Otherwise, $\ast\alpha\in t(w_0)$. For all $n$, $(w_n,w_{n+1})$ is sensible so $\alpha\in t(w_n)$ and 
$\<\lim \mathfrak D\upharpoonright Y,{\sigma}^n\({\vec w}\)\>\models\alpha;$
hence
$\<\lim \mathfrak D\upharpoonright Y,\({\vec w}\)\>\models\ast\alpha.$
\endproof

We are now ready to prove the main theorem of this section, which in particular implies that our semantics are sound for $\mathcal{DTL}$.

\begin{thm}\label{second}
Let $\varphi$ be a formula of $\mathcal{L}$, and suppose $\varphi$ is satisfied in a non-deterministic quasimodel
$\mathfrak{ D}=\<W,R,g,t\>.$
Then, there exists ${{{\vec w}^*}\in W^g}$ such that
$\<\lim\mathfrak{D},{{\vec w}}^*\>\models\varphi.$
\end{thm}
\proof
Pick $w^*\in W$ such that $\varphi\in t(w^*)$. By Lemma \ref{extension}, $w^*$ can be included in a realizing path ${\vec w}^*$. It follows from Lemma \ref{sound} that
$\<\lim\mathfrak{D},{\vec w}^*\>\models\varphi.$
\endproof
\section{Local Kripke frames}\label{secfour}

In this section we establish a basic framework for describing small substructures of Kripke frames. We wish to work with locally finite frames, and often it is convenient to give explicit bounds on the size of neighborhoods. These bounds will depend on the length of $\varphi$, denoted ${\left|\varphi\right|}.$

\begin{defn}[local Kripke frame]
A {\em local $\varphi$-typed Kripke frame} is a tuple
${\mathfrak {a}}=\<w_{\mathfrak {a}},W_{\mathfrak {a}},R_{\mathfrak {a}},t_{\mathfrak {a}}\>,$
where $\<W_{\mathfrak {a}},R_{\mathfrak {a}},t_{\mathfrak {a}}\>$ is a $\varphi$-typed Kripke frame and $w_{\mathfrak {a}}\in W_{\mathfrak {a}}$ is such that $R_{\mathfrak a}w_{\mathfrak {a}} v$ for all $v\in W_{\mathfrak {a}}$.
\end{defn}

The reader may recognize local Kripke frames as being nothing more than Kripke frames with a root.  The reason we call them `local' here is that, for our purposes, $\mathfrak {a}$ will represent a neighborhood of $w_{\mathfrak {a}}$, which may be a world in a larger Kripke frame $\mathfrak{A}$. The frame $\mathfrak{A}$ will often be disconnected and, hence, have no candidate for a root. The lower-case letters used to denote local Kripke frames are meant to be suggestive of this local character.

We will write $t(\mathfrak a)$ instead of $t_{\mathfrak {a}}\(w_{\mathfrak {a}}\)$.

\subsection{Tree-like Kripke frames}\label{secfive}

Given a local Kripke frame
${\mathfrak {a}}=\<w,W,R,t\>,$
the relation $R$ induces an equivalence relation $\sim_R$ on $W$ given by
$w\sim_R v\Leftrightarrow Rwv\text{ and }Rvw.$

The equivalence class of a world $w$ is usually called the {\em cluster} of $w$. We will denote it by $[w]_R$, or simply $[w]$ if this does not lead to confusion. $R$ then induces a partial order on $W / \sim_R$ defined by
$R[w][v]\Leftrightarrow Rwv.$
If this partial order forms a tree (that is, if whenever $R[u][w]$ and $R[v][w]$, then either $R[u][v]$ or $R[v][u]$), we will say that ${\mathfrak {a}}$ is {\em tree-like}.

Given a tree-like local Kripke frame $\mathfrak a$, we can define $\mathrm{hgt}(\mathfrak a)$ and $\textrm{wdt}(\mathfrak a)$ as the height and width of $\mathfrak a/\sim_R$. Likewise, we will define the {\em depth} of $\mathfrak a$, $\mathrm{dpt}(\mathfrak a)$, to be the maximum number of elements in a single cluster of $W_{\mathfrak a}$.

\begin{defn}[norm of a Kripke frame]
Let $\mathfrak a$ be a tree-like local Kripke frame. Define the {\em norm} of $\mathfrak a$, denoted $\left\|\mathfrak a\right\|$, by
\[\left\|\mathfrak a\right\|=\max(\mathrm{hgt}(\mathfrak a),\mathrm{wdt}(\mathfrak a),\mathrm{dpt}(\mathfrak a)).\]
\end{defn}

We will use the norm of a local Kripke frame as a measure of its size rather than the more obvious $\left|W_{\mathfrak a}\right|$, because it is often more manageable. However, it is clear that one can use the norm of a frame to find bounds for the number of worlds in it (and vice-versa).

For the rest of this paper, all local Kripke frames will be assumed to be tree-like.

\subsection{Binary relations between local Kripke frames}

Many times it will be useful to compare different local Kripke frames and express relations between them. The following binary relations are essential and will appear throughtout the text:

\begin{defn}[reduction of local Kripke frames]
Say that ${\mathfrak {b}}$ {\em reduces} to ${\mathfrak {a}}$ (or, alternately, $\mathfrak {a}$ {\em embeds} into $\mathfrak {b}$), denoted ${\mathfrak {a}}\unlhd{\mathfrak {b}}$, if there exists an injective function
$e:W_{\mathfrak {a}}\to W_{\mathfrak {b}}$
such that, for all $w,v\in W_{\mathfrak {a}}$,
$R_{\mathfrak {a}} wv\Leftrightarrow R_{\mathfrak {b}} e(w)e(v),$
$t_{\mathfrak {a}}(w)=t_{\mathfrak {b}}(e(w))$ and $e\(w_{\mathfrak {a}}\)=w_{\mathfrak {b}}$.
\end{defn}

Roughly, if ${\mathfrak {a}}\lhd{\mathfrak {b}}$, $\mathfrak {b}$ contains worlds which could be removed without altering $t(\mathfrak b)$, and hence $\mathfrak {b}$ could be replaced by $\mathfrak {a}$ for most purposes. We will make this precise later in this section.

\begin{defn}[subframe]
For $v\in W_{\mathfrak {a}}$, set
${\mathfrak {a}}^v=\<v,W^v_{\mathfrak {a}},R^v_{\mathfrak {a}},t^v_{\mathfrak {a}}\>,$
where
$W^v_{\mathfrak {a}}=\cbra w\in W_{\mathfrak {a}}:R_{\mathfrak {a}} vw \cket$
and $R^v_{\mathfrak {a}},t^v_{\mathfrak {a}}$ are the corresponding restrictions to $W^v_{\mathfrak {a}}$.

Then define
${\mathfrak {b}}\preceq{\mathfrak {a}}$
if ${\mathfrak {b}}={\mathfrak {a}}^v$ for some $v\in W_{\mathfrak {a}};$ we will say ${\mathfrak {b}}$ is a {\em subframe} of ${\mathfrak {a}}$.
\end{defn} 

If ${\mathfrak {a}}\preceq{\mathfrak {b}}$ and ${\mathfrak {b}}\preceq{\mathfrak {a}}$, we will write ${\mathfrak {a}}\sim{\mathfrak {b}}$. Likewise, ${\mathfrak {a}}\prec{\mathfrak {b}}$ means that ${\mathfrak {a}}\preceq{\mathfrak {b}}$ but not vice-versa, while ${\mathfrak {a}}\prec_1{\mathfrak {b}}$ means that ${\mathfrak {a}}\prec{\mathfrak {b}}$ and there is no intermediate local Kripke frame ${\mathfrak{c}}$ such that ${\mathfrak {a}}\prec{\mathfrak{c}}\prec{\mathfrak {b}}$.

Suppose ${\mathfrak {b}}\preceq{\mathfrak {a}}$. If we think of ${\mathfrak {a}}$ as a neighborhood of $w_{\mathfrak {a}}$, then ${\mathfrak {b}}$ represents an open subset of $W_{\mathfrak {a}}$ (which does not necessarily contain $w_{\mathfrak {a}}$). If it does contain $w_{\mathfrak {a}}$, then $\mathfrak {b}\sim\mathfrak {a}$, and the two represent the same open set but maybe `centered' at a different point.

\begin{defn}[subframe representatives]
Let $\mathfrak a$ be a local Kripke frame. A {\em set of subframe representatives for} $\mathfrak a$ is a set of representatives of the equivalence classes of
$\cbra\mathfrak b:\mathfrak b\prec_1\mathfrak a\cket$
under $\sim$.
\end{defn}

In the following definition and throughout the paper, a binary relation $g$ between Kripke frames is {\em non-confluent} if whenever $gwv$, $gw'v'$ and $Rvv'$, it follows that $Rww'$, so that we can fill in the dotted arrow in the following diagram:
$$
\xymatrix
{
w'\ar[r]^g&v'\\
w\ar@{-->}[u]^R\ar[r]^g&v\ar[u]^R.
}
$$

\begin{defn}[termporal successor]\label{ts}
Say $\mathfrak a$ is an {\em temporal successor} of $\mathfrak b$, denoted $\mathfrak a\rra\mathfrak b$, if there exists a non-confluent sensible relation
$g\subseteq W_{\mathfrak a}\times W_{\mathfrak a}$
such that $gw_{\mathfrak a}w_{\mathfrak b}$.
\end{defn}

\begin{lem}\label{reduction}
Let $\varphi$ be any formula of $\mathcal{DTL}$, and $\mathfrak a\rra\mathfrak b$ be local Kripke frames.

Then, there exists $\mathfrak d\unlhd\mathfrak b$ such that $\mathfrak a\rra\mathfrak d$ and
$\|\mathfrak d\|\leq\|\mathfrak a\|+\left|\varphi\right|.$
\end{lem}
\proof
We will skip the proof. The general idea is that if $\|\mathfrak b\|>\|\mathfrak a\|+|\varphi|$, then $W_{\mathfrak b}$ contains worlds which could be deleted, giving us $\mathfrak d$ such that
$\mathfrak a\rra\mathfrak d\lhd\mathfrak b.$ Repeating this enough times we can attain the desired bound.
\endproof

\subsection{The space of bounded frames}\label{ikf}

\begin{defn}[$\mathfrak I_K(\varphi)$]
Let $\varphi$ be any formula of $\mathcal{DTL}$ and $K\geq 0$. Define $I_K(\varphi)$ to be the set of all local, tree-like Kripke frames $\mathfrak a$ such that
$\|\mathfrak a\|\leq (K+1)|\varphi|.$
\end{defn}

Now, consider
$I_\omega(\varphi)=\bigcup_{k\geq 0}I_k(\varphi)$
(evidently this is the set of all finite local tree-like frames).
Define
$\mathfrak{I}_\omega(\varphi)=\<I_\omega(\varphi),\succeq,\rra,t\>,$
where $t(\mathfrak a)={\bf t}_{\mathfrak a}$.

$\mathfrak I_\omega(\varphi)$ is a $\varphi$-typed Kripke frame\footnote
{
Note that the accessibility relation is written as $\succeq$, so that $\Box\psi$ holds in $\mathfrak a$ if $\psi$ holds in $\mathfrak b$ for all $\mathfrak b\preceq\mathfrak a$, even though it is more standard to write the accessibility relation in the opposite direction. We believe it is natural to adopt this convention because $\|\mathfrak b\|\leq\|\mathfrak a\|$ whenever $\mathfrak b\preceq\mathfrak a$.
}
with a sensible relation $\rra$; however, it is not necessarily $\omega$-sensible, so $\mathfrak{I}_\omega(\varphi)$ is not a non-deterministic quasimodel as it stands. It does contain substructures which are, as we will see in the next section.

\subsection{Building local Kripke frames from subframes}

Often we will want to construct a local Kripke frame from smaller pieces. Here we will define the basic operation we will use to do this, and establish the conditions that the pieces must satisfy.
\begin{defn}
Let $\varphi$ be a formula of $\mathcal{DTL}$, $T\subseteq{\bf type}(\varphi)$ and $A\subseteq I_\omega(\varphi)$.

For each ${\bf t}\in T$ define
$\bra\,T\oplus A\,\ket_{\bf t}=\<w,W,R,t\>$
by setting $w={\bf t},$
\begin{align*}
W&=T\cup\coprod_{\mathfrak a\in A}W_{\mathfrak a},\\
R&=\(T\times W\)\cup\coprod_{\mathfrak a\in A}R_{\mathfrak a}
\end{align*}
and
\[
t(w)=
\begin{cases}
w&\text{ if $w\in T$}\\
t_{\mathfrak a}(w)&\text{ if $w\in W_{\mathfrak a}$}.
\end{cases}
\]
\end{defn}

\begin{figure}[htp]
\begin{center}
\scalebox
{0.9}
{
$$
\xymatrix{
\mathfrak a_0&\mathfrak a_1&\mathfrak a_2&...&\mathfrak a_J\\
&&\cbra {\bf t}_m\cket_{m\leq M}\ar@{~>}[ul]_R\ar@{~>}[ull]_R\ar@{~>}[u]_R\ar@{~>}[urr]_R&&
}
$$
}
\end{center}
\caption
{
If $A=\cbra\mathfrak a_j\cket_{j\leq J}$ and $T=\cbra {\bf t}_m\cket_{m\leq M}$, $[\,T\oplus A\,]_{{\bf t}_0}$ has ${\bf t}_0$ as a root and each $\mathfrak a_j$ as a subframe.
}
\end{figure}
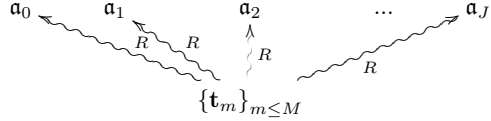
\begin{defn}[admitting pair]\label{admit}
Let $\varphi$ be a formula of $\mathcal{DTL}$, $T\subseteq {\bf type}(\varphi)$, ${\bf t}\in T$, $\mathfrak{c}\in I_\omega(\varphi)$ and $A\subseteq I_\omega(\varphi)$.

The triple $\<\,T,A,{\bf t}\,\>$ {\em admits} $\mathfrak b$ if $(t(\mathfrak b),{\bf t})$ is sensible and either
\begin{enumerate}
\item there is $\mathfrak a\in A$ such that $\mathfrak b\rra\mathfrak a$ or
\item
\begin{enumerate}
\item for each $w\in[w_{\mathfrak b}]$ there exists ${\bf s}\in T$ such that $(t(w),{\bf s})$ is sensible and
\item there is a set of subframe representatives $B$ for $\mathfrak b$ and an injection
$\iota:B\to A$
such that for each $\mathfrak c\in B$, $\mathfrak c\rra\iota(\mathfrak c)$.
\end{enumerate}
\end{enumerate}
\end{defn}
\begin{lem}\label{tempsucc}
If $\mathfrak a$ and $\mathfrak b$ are local Kripke frames, $B$ is a set of subframe representatives for $\mathfrak b$ and the triple
$\<\,t_{\mathfrak b}[w_{\mathfrak b}],B,t(\mathfrak b)\,\>$
admits $\mathfrak a$, then $\mathfrak a\rra\mathfrak b$.
\end{lem}
\proof
The proof is straightforward and we omit it here.
\endproof

\begin{defn}[coherence]\label{cohe}
Let $\varphi$ be a formula of $\mathcal{DTL}$, $T\subseteq{\bf type}(\varphi)$, ${\bf t}\in T$ and $A\subseteq I_\omega(\varphi)$.

Then,
\begin{enumerate}
\item the pair $\<\,T,A\,\>$ is {\em coherent} if, for all $\Box\psi\in {\mathrm sub}_\pm(\varphi)$ and ${\bf t}\in T$, $\Box\psi\in {\bf t}$ if and only if $\psi\in{\bf s}$ for all ${\bf s}\in T$ and $\Box\psi\in t(\mathfrak a)$ for all $\mathfrak a\in A$ and
\item if $\mathfrak c\in I_\omega(\varphi)$, the pair $\<\,T,A\,\>$ is {\em coherent for $\mathfrak c$} if it is coherent and admits $\mathfrak c$.
\end{enumerate}
\end{defn}

Coherent pairs are useful for constructing local Kripke frames.

\begin{lem}\label{add}
Let $\varphi$ be a formula of $\mathcal{DTL}$, $T\subseteq{\bf type}(\varphi)$, $A\subseteq I_\omega(\varphi)$ and ${\bf t}\in T$.

Then,
\begin{enumerate}
\item $\bra\,T\oplus A\,\ket_{\bf t}$ is a local Kripke frame if and only if $\<\,T,A\,\>$ is coherent and
\item $\bra\,T\oplus A\,\ket_{\bf t}$ is a local Kripke frame and $\mathfrak c\rra\mathfrak a$ whenever $\<\,T,A\,\>$ is coherent for $\mathfrak c$.
\end{enumerate}
\end{lem}
\proof
To prove 1, one can check that the coherence conditions correspond exactly to the condition in Definition \ref{frame}. The second claim follows from the first and Lemma \ref{tempsucc}.
\endproof
\section{Simulating topological models}\label{secsim}
In this section we will study simulations, which are the basic tool for extracting non-deterministic quasimodels from dynamic topological models.

If
$\mathfrak{M}=\<X,\mathcal T,f,V\>$
is a dynamic topological model and $x\in X$, assign a $\varphi$-type $\tau(x)$ to $x$ given by
$\tau(x)=\cbra\psi\in \mathrm{sub}_\pm(\varphi):x\in V(\psi)\cket.$
We will also define
$\tau^\Diamond(x)=\cbra\psi\in \mathrm{sub}_\pm(\varphi):\Diamond\psi\in \tau(x)\cket.$

Analogously, if
$\mathfrak F=\<W,R,t\>$
is a $\varphi$-typed Kripke frame and $w\in W$, set
$t^\Diamond(w)=\cbra\psi\in \mathrm{sub}_\pm(\varphi):\Diamond\psi\in t(w)\cket.$
\subsection{Simulations}
\begin{defn}[simulation]
Let $\varphi$ be a formula of $\mathcal{DTL}$,
\[\mathfrak{M}=\<X,\mathcal T,f,V\>\]
a dynamic topological model and
$\mathfrak F=\<W,R,t\>$
a $\varphi$-typed Kripke frame.

A continuous relation
$\chi\subseteq W\times X$
is a {\em simulation} if, for all $x\in X,$
\[x\in\chi(w)\Rightarrow \tau(x)=t(w).\]
\end{defn}

We will call the latter property {\em type-preservation}.

A simulation can be thought of as a one-way bisimulation; a topological bisimulation would be an open, continuous map which preserves valuations of propositional variables. Simulations can be used to capture much of the purely topological information about $\mathfrak{M}$. However, temporal behavior is disregarded here; for this we need simulations on non-deterministic quasimodels, not just Kripke frames, and these simulations must respect the transition function.

\begin{defn}[$\omega$-simulation]\label{simulation}
Let $\varphi$ be a formula of $\mathcal{DTL}$ and
$\mathfrak{M}=\<X,\mathcal T,f,V\>$
a dynamic topological model. Let
$\mathfrak{F}=\<W,R,t\>$
be a Kripke frame and $g$ a sensible relation on $W$.

Suppose
$\chi\subseteq W\times X$
is a simulation.

Then, $\chi$ is an $\omega$-simulation if
$
f\chi\subseteq\chi g.
$
\end{defn}

While $g$ is not required to be $\omega$-sensible on $\mathfrak{F}$, we can use $\chi$ to extract a non-deterministic quasimodel from $\mathfrak{F}$.

\begin{lem}\label{isquasi}
Suppose $\mathfrak{F}=\<W,R,t\>$ is a $\varphi$-typed Kripke frame with a sensible relation $g$, $\mathfrak{M}=\<X,\mathcal T,f,V\>$ is a dynamic topological model and $\chi\subseteq W\times X$ is an $\omega$-simulation.

Then, $\mathfrak F\upharpoonright \mathrm {dom}(\chi)$ is a non-deterministic quasimodel.
\end{lem}
\proof We only need to prove that $g\upharpoonright \mathrm {dom}(\chi)$ is $\omega$-sensible.

Let $w\in \mathrm {dom}(\chi)$, $\neg\ast\psi\in t(w)$ and $x\in \chi(w)$.

Then,
$\<\mathfrak M,x\>\models\neg\ast\psi,$ so for some $N>0$,
$f^{N}(x)\models\neg\psi;$
but $f\chi\subseteq\chi g,$ so there exists $v\in W$ such that $f^N(x)\in \chi(v)$ (hence $v\in \mathrm {dom}(\chi)$) and $g^Nwv$.

Thus $\neg\psi\in t(v)$, which is what we wanted.
\endproof

Suppose that $\chi$ is an $\omega$-simulation and $x^\ast \in X$ is such that
$\<\mathfrak M,x^\ast\>\models\varphi.$
If there exists $w^\ast\in W$ such that $x^\ast\in\chi(w^\ast)$, then clearly $\mathfrak{D}$ satisfies $\varphi$, since $\varphi\in t(w^\ast)$.

Thus, if we show that, given a dynamic topological model
$\mathfrak{M}=\<X,\mathcal T,f,V\>,$
there exists a non-deterministic quasimodel
$\mathfrak{D}=\<W,R,g,t\>$
with a surjective $\omega$-simulation
$\chi\subseteq W\times X,$
this would imply that, given any satisfiable formula, it can be satisfied in a non-deterministic quasimodel.

In fact, we will show that $\mathfrak I_\omega(\varphi)$ (defined in Section \ref{ikf}) contains a `canonical' quasimodel in the sense that there always exists a surjective $\omega$-simulation from $I_\omega(\varphi)$ to $X$.
\begin{lem}
Given any dynamical topological model
$\mathfrak{M}=\<X,\mathcal{T},f,V\>,$
there exists a unique maximal simulation
$\chi^\ast\subseteq I_\omega(\varphi)\times X.$
\end{lem}
\proof
The proof uses Zorn's Lemma in a straightforward fashion and we skip it. For uniqueness, if $\chi^\ast$ and $\chi^+$ are two maximal simulations, one can check that $\chi^\ast\cup\chi^+$ is also a simulation, which must equal both $\chi^\ast$ and $\chi^+$ by maximality.
\endproof

Given a simulation
$\chi\subseteq I_\omega(\varphi)\times X,$
we will denote $\chi\upharpoonright{I_k(\varphi)}$ by $\chi_k$.

Our goal is to prove that $\chi^\ast$ gives us a surjective $\omega$-simulation. The following lemma will be essential in proving this.
\begin{lem}[simulation extensions]\label{simex}
Let
$\mathfrak M=\<X,\mathcal T,f,V\>$
be a dynamic topological model and
$\chi\subseteq I_\omega(\varphi)\times X$
be a simulation.

If $T\subseteq{\bf type}(\varphi)$ and $A\subseteq I_\omega(\varphi)$ are coherent (as in Definition \ref{cohe}) and there is a set $E\subseteq X$ such that, for all $\mathfrak a\in A$,
$E\subseteq \overline{\chi(\mathfrak a)}$
and, for all ${\bf t}\in T$,
$E\subseteq \overline{E\cap V\(\bigwedge {\bf t}\)},$
then
\[\zeta=\chi\cup\cbra\left ([\,T\oplus A\,]_{\tau(x)},x \right ):x\in E\text{ and }\tau(x)\in T\cket\]
is also a simulation.
\end{lem}
\proof
Note, first, that $\zeta$ is type-preserving. We must prove it is also continuous.

Consider an arbitrary open set $U\subseteq X$. If $U\cap E=\varnothing$,
$\zeta^{\,-1}(U)=\chi^{-1}(U),$
which is open because $\chi$ is continuous.

Otherwise,
$\zeta^{\,-1}(U)=\chi^{-1}(U)\cup \zeta^{\,-1}(E\cap U).$

We must prove that if $\mathfrak b\in \zeta^{\,-1}(E\cap U)$ and $\mathfrak c\preceq\mathfrak b$, then
$\mathfrak c\in \chi^{-1}(U)\cup \zeta^{\,-1}(E\cap U).$
First assume that $\mathfrak c\prec \mathfrak b$. In this case, $\mathfrak c\preceq\mathfrak a$ for some $\mathfrak a\in A$. 

Since
\[E\subseteq \overline{\chi(\mathfrak a)},\]
there exists some $y\in U\cap\chi(\mathfrak a)$. Because $U$ is open and $\chi$ is continuous, there also exists $z\in U\cap\chi(\mathfrak c)$, as we wanted.

Otherwise, $\mathfrak c\sim\mathfrak b$. Now,
$E\subseteq \overline{E\cap V\(\bigwedge t(\mathfrak c)\)},$
so there exists $y\in E\cap U$ such that $\tau(y)=t(\mathfrak c)$. By the definition of $\zeta$, $y\in \chi(\mathfrak c)$, so $\mathfrak c\in\zeta^{\,-1}(U)$.

We conclude that $\zeta$ is continuous, and therefore a simulation.
\endproof

\begin{prop}\label{total}
For any dynamic topological model $\mathfrak{M}$, $\chi^\ast_0$ is surjective.
\end{prop}
\proof
We must define slightly stronger bounds on frames for our proof to go through; namely, define {\em small} and {\em very small} as follows:
\begin{enumerate}
\item Say $\mathfrak a$ is {\em very small} if
$\|\mathfrak a\|\leq \left|t^\Diamond(\mathfrak a)\right|$
and $\mathrm{hgt}(\mathfrak a)<\left|t^\Diamond(\mathfrak a)\right|.$
\item Say $\mathfrak a$ is {\em small} if
$\|\mathfrak a\|\leq \left|t^\Diamond(\mathfrak a)\right|$
and there is a very small $\mathfrak b\preceq\mathfrak a$ such that
$t^\Diamond(\mathfrak b)=t^\Diamond(\mathfrak a).$
\end{enumerate}

Note that if $\mathfrak a$ is small, then $\mathfrak a\in I_0(\varphi)$. We will prove that $\chi^\ast$ is surjective even when restricted to the set of small frames.

Suppose
$\chi\subseteq I_\omega(\varphi)\times X$
is any simulation.

Say a point $x\in X$ is {\em bad} if there is no small $\mathfrak a$ such that $x\in\chi(\mathfrak a)$.

We claim that $\chi$ is not maximal if there are bad points.

To see this, let $E$ be an arbitrary subset of $X$. Define
\[\mathrm{Bad}(E)=\cbra \tau(x):x\in E\text{ is bad}\cket.\]

Assume that
$\mathrm{Bad}(X)\not=\varnothing.$

Then, we can pick out an open set $U_*$ which minimizes
$\left|\mathrm{Bad}(U)\right|+\left|\chi^{-1}_0(U)\right|,$
where $U$ ranges over all open sets that contain bad points. Notice that $\mathrm{Bad}(E)$ and $\chi^{-1}_0(E)$ are finite, since they are subsets of ${\bf type}(\varphi)$ and $I_0(\varphi)$, which are finite.

Note also that, for such a $U_\ast$, whenever $U\subseteq U_\ast$ is open and contains bad points, we know that $\mathrm{Bad}(U)=\mathrm{Bad}(U_\ast)$ and $\chi^{-1}_0(U)=\chi^{-1}_0(U_\ast)$ (otherwise $U_\ast$ would not be optimal).

We will construct a simulation $\zeta\supseteq\chi$ such that $\zeta_0\not=\chi_0$.

Let $\Psi$ be the set of all formulas $\psi\in\mathrm{sub_{\pm}}(\varphi)$ such that
$\Diamond\psi\in \bigcup\mathrm{Bad}(U_\ast)$
but
$\psi\not\in \bigcup\mathrm{Bad}(U_\ast).$

For each $\psi\in \Psi$, $U$ contains a point $y$ such that
$\<\mathfrak M,y\>\models\psi.$
Note that $y$ cannot be bad, so there exists a small frame $\mathfrak c$ such that $y\in\chi_0\(\mathfrak c\),$ and hence a very small $\mathfrak a\preceq\mathfrak c$. Set $\mathfrak a=\mathfrak a_\psi$, and
$A=\cbra \mathfrak a_\psi:\psi\in\Psi\cket.$

Pick a minimal, non-empty $T\subseteq \mathrm{Bad}(U_\ast)$ such that $T$ and $A$ are coherent. By Lemma \ref{add},
$\mathfrak a_\ast=[\,T\oplus A\,]_{{\bf t}_\ast}$
is local Kripke frame, and we can set
\[\zeta=\chi\cup \cbra (\mathfrak b,y):y\in U_\ast,\mathfrak b\sim\mathfrak a\text{ and }\tau(y)=t(\mathfrak b)\cket.\]

By Lemma \ref{simex}, $\zeta$ is a simulation. It remains to show that $\mathfrak a_\ast$ is small.

Note first that, since all elements of $A$ are very small, for all $\mathfrak a\in A$,
$\mathrm{hgt}(\mathfrak a)<\left|t^\Diamond(\mathfrak a)\right|.$

This shows that
$\mathrm{hgt}(\mathfrak a_\ast)\leq\left|t^\Diamond(\mathfrak a_\ast)\right|,$
and equality holds only if
$t^\Diamond(\mathfrak b)=t^\Diamond(\mathfrak a_\ast)$
for some $\mathfrak b\in A$; this gives us the very small frame $\mathfrak b\preceq \mathfrak a_\ast$.

Similarly,
$\mathrm{wdt}(\mathfrak a)\leq\left|t^\Diamond(\mathfrak a_\ast)\right|$
for all $\mathfrak a\in A$, and $\mathfrak a_\ast$ has at most $\left|t^\Diamond(\mathfrak a_\ast)\right|$ immediate successors, so
$\mathrm{wdt}(\mathfrak a_\ast)\leq\left|t^\Diamond(\mathfrak a_\ast)\right|.$
One can easily see that
$\mathrm{dpt}(\mathfrak a_\ast)\leq|\varphi|.$

It follows that $\mathfrak a_\ast\in I_0(\varphi)$. Since $U_\ast$ contained bad points, $\chi\subsetneq \zeta$, as desired.
\endproof
\begin{prop}\label{temptotal}
For all $K\geq 0$,
$f\chi^\ast_K\subseteq\chi^\ast_{K+1}g.$
\end{prop}
Note that, as a consequence of this,
$f\chi^\ast\subseteq\chi^\ast g$
and thus $\chi^\ast$ is an $\omega$-simulation.
\proof
The proof follows much the same structure as that of Proposition \ref{total}.

Suppose $\mathfrak a\rra \mathfrak b$. Define {\em small relative to $\mathfrak a$} and {\em very small relative to $\mathfrak a$} as follows:
\begin{enumerate}
\item Say $\mathfrak b$ is {\em very small relative to }$\mathfrak a$ if
$\|\mathfrak b\|\leq \|\mathfrak a\|+|\varphi|$
and $\mathrm{hgt}(\mathfrak a)<\|\mathfrak a\|+|\varphi|.$
\item Say $\mathfrak b$ is {\em small relative to }$\mathfrak a$ if
$\|\mathfrak b\|\leq \|\mathfrak a\|+|\varphi|$
and there is $\mathfrak d\preceq\mathfrak b$ which is very small relative to $\mathfrak a$.
\end{enumerate}

Let $\chi$ be a simulation.

Say $x\in X$ {\em fails} for $\mathfrak a$ if
$f^{-1}(x)\cap\chi(\mathfrak a)\not=\varnothing,$
but there is no $\mathfrak b$ which is small relative to $\mathfrak a$ such that $x\in\chi(\mathfrak b)$.

We claim that if $\chi$ contains points that fail for any $\mathfrak a$, then $\chi$ is not maximal.

Suppose there exists $\mathfrak a_\ast$ such that some point fails for $\mathfrak a_\ast$, and $\mathfrak a_\ast$ is minimal with this property.

For $E\subseteq X$ and $\mathfrak a\in I_\omega(\varphi)$ define
\[\mathrm{Fail}(E)=\cbra\tau(x):x\in E\text{ fails for some }\mathfrak c\sim\mathfrak a_\ast\cket.\]

Pick an open set $U_\ast$
which minimizes
$\left|\mathrm{Fail}(U)\right|+\left|\chi_{K+1}^{-1}(U)\right|,$
where $U$ ranges over all open sets which contain points that fail for $\mathfrak a_\ast$. 

As before, let $\Psi$ be the set of all formulas $\psi\in\mathrm{sub_{\pm}}(\varphi)$ such that
$\Diamond\psi\in \bigcup\mathrm{Fail}(U_\ast)$
but
$\psi\not\in \bigcup\mathrm{Fail}(U_\ast).$
To each element $\psi$ of $\Psi$ assign a very small frame
$\mathfrak b_\psi\in \chi^{-1}(U_\ast)$
such that $\psi\in t\(\mathfrak b_\psi\)$. These frames exist by Lemma \ref{total}.

Pick any ${\bf t}_\ast\in \mathrm{Fail}(U_\ast)$ such that $(t(\mathfrak a_\ast),{\bf t}_\ast)$ is sensible.

We must find $T$ and $B'$ such that the triple $\<T,B',{\bf t}_\ast\>$ admits $\mathfrak a_\ast$. Here we will consider two cases.

\noindent{\bf Case 1.} Suppose there is $\mathfrak c\sim\mathfrak a_\ast$ such that no point of $U_\ast$ fails for $\mathfrak c$. In this case, there must exist
$\mathfrak d_\ast\in\chi_{K+1}^{-1}(U)$ which is very small relative to $\mathfrak c$.

If this holds, set $B'=\cbra\mathfrak d_\ast\cket$ and $T=\cbra {\bf t}_\ast\cket$.

\noindent{\bf Case 2.} Suppose Case 1 does not hold. Note that in this case, given any $\mathfrak c\sim\mathfrak a_\ast$, there is some point $z$ in $U_\ast$ which fails for $\mathfrak c$, and hence there is ${\bf t}_{\mathfrak c}=\tau(z)\in \mathrm{Fail}(U_\ast)$ such that $(t(\mathfrak c),{\bf t}_{\mathfrak c})$ is sensible. Set $T=\cbra {\bf t}_{\mathfrak c}:\mathfrak c\sim\mathfrak a_\ast\cket$.

Let $C$ be a set of subframe representatives for $\mathfrak a_\ast$. For each ${\mathfrak{c}}\in C$, no point of $U_\ast$ fails for ${\mathfrak{c}}$ (because we picked $\mathfrak a_\ast$ to be minimal). However, $f^{-1}(U_\ast)$ is open and therefore contains points in $\chi(\mathfrak c)$, so there must exist a frame $\mathfrak b_{\mathfrak c}$ which is small relative to $\mathfrak c$ (and, by passing to a subframe if necessary, we can pick it to be very small relative to $\mathfrak c$).

Then define
$B'=\cbra \mathfrak b_{\mathfrak c}:\mathfrak c\in C\cket.$

In either of the two cases set
$B=B'\cup\cbra \mathfrak b_\psi\cket_{\psi\in\Psi}.$

Then, by Lemma \ref{add}, \[\mathfrak b_\ast=[\,T\oplus B\,]_{{\bf t}_\ast}\]
is a local Kripke frame and $\mathfrak a_\ast\rra\mathfrak b_\ast$ by Lemma \ref{tempsucc}.

We can then set
\[\zeta=\chi\cup \cbra (\mathfrak d,y):y\in U_\ast,\mathfrak d\sim\mathfrak b_\ast\text{ and }\tau(y)=t(\mathfrak d)\cket.\]

One can then show as before that $\mathfrak b_\ast$ is small for $\mathfrak a_\ast$, so $\zeta$ is a simulation which properly contains $\chi$. Therefore, $\chi$ is not maximal.
\endproof

We are now ready to give a completeness proof of non-deterministic semantics for $\mathcal{DTL}$.

\begin{defn}
Given a dynamic topological model $\mathfrak M$ satisfying $\varphi$, define
$\mathfrak M/\varphi=\mathfrak I_\omega(\varphi)\upharpoonright \mathrm{dom}(\chi^\ast).$
\end{defn}

\begin{thm}
If $\mathfrak M$ is a dynamic topological model satisfying $\varphi$, then $\mathfrak M/\varphi$ is a non-deterministic quasimodel satisfying $\varphi$.
\end{thm}
\proof
By Proposition \ref{temptotal}, $\chi^\ast$ is an $\omega$-simulation, so by Lemma \ref{isquasi}, ${\mathfrak M}/\varphi$ is a $\varphi$-typed non-deterministic quasimodel.

Pick $x_\ast\in X$ such that
$\<\mathfrak M,x_\ast\>\models \varphi.$
By Proposition \ref{total}, $\chi^\ast$ is surjective, so there exists $\mathfrak a_\ast\in I_\omega(\varphi)$ such that $x_\ast\in\chi^\ast(\mathfrak a_\ast)$; hence $\varphi\in t(\mathfrak a_\ast)$.

This shows that that $\mathfrak M/\varphi$ satisfies $\varphi$, as desired.
\endproof
\section{A model-search procedure}\label{re}
Non-deterministic quasimodels can be used to give a recursive enumeration of all valid formulas of $\mathcal{DTL}$. The general strategy is to generate finite `chunks' of $\varphi$-typed non-deterministic quasimodels; if the search for chunks of arbitrary size terminates, $\varphi$ is not satisfiable. Otherwise we can construct a non-deterministic quasimodel for $\varphi$ and hence a model.

We will use Kruskal's Tree Theorem to give a recursive enumeration of all valid formulas of $\mathcal{DTL}$. Most of what follows is an adaptation of a proof in \cite{konev} that $\mathcal{DTL}_1$, the fragment of $\mathcal{DTL}$ where $\ast$ is not allowed to appear in the scope of $\Box$, is recursively enumerable. We will use non-deterministic quasimodels to generalize this result to full $\mathcal{DTL}$.

Recall that a pair $\<S,\leq\>$ is a well-partial order if, for any infinite sequence
$\< s_n\>_{n\geq 0}\subseteq{S},$
there exist indices $M_0<M_1$ such that $s_{M_0}\leq s_{M_1}$.

A labeled tree is a triple $\<T,\leq,L\>$, where $\<T,\leq\>$ is a tree and $L:T\to \Lambda$ is a labeling function to some fixed set $\Lambda$ of labels.

If $\mathfrak T_0$ and $\mathfrak T_1$ are labeled trees and $\Lambda$ is partially ordered, an embedding between $\mathfrak T_0$ and $\mathfrak T_1$ is a function $e:{ T}_0\to { T}_1$ which is an embedding as trees and such that $L_0\leq_\Lambda L_1e$. If such an embedding exists, we say that $\mathfrak T_0\leq\mathfrak T_1$. We will always assume that embeddings map roots to roots.

\begin{thm}[Kruskal]
The set of finite trees with labels in a well partially-ordered set is well-partially ordered.
\end{thm}
\proof
Kruskal's original proof can be found in \cite{ktt}.
\endproof
We wish to apply Kruskal's Theorem to elements of $I_{\omega}(\varphi)$.
\begin{lem}\label{kruskapp}
Let
${\mathfrak {a}}_0,{\mathfrak {a}}_1,...,{\mathfrak {a}}_n,...$
be an infinite sequence of finite, tree-like local Kripke frames.

Then, there exist $M_1<M_2$ such that ${\mathfrak {a}}_{M_1}\unlhd{\mathfrak {a}}_{M_2}$.
\end{lem}
\proof
This is a straightforward application of Kruskal's Tree Theorem; we only need to represent tree-like local Kripke frames by labeled trees.

Namely, to each tree-like local Kripke frame
${\mathfrak {a}}=\<w,W,R,t\>$
assign a labeled tree
$\mathfrak T_{\mathfrak a}=\<T_{\mathfrak a},\leq_{\mathfrak {a}},l_{\mathfrak {a}}\>$
given by $\mathfrak T_{\mathfrak {a}}=W/\sim_R$ and $\leq_{\mathfrak {a}}=R/\sim_R$. If $[w]=\cbra w_i\cket_{i\leq I}$, let
$l([w])=\<t(w_0),...,t(w_I)\>.$

It is known that the set of finite sequences of elements of a finite set is well-partially ordered by the `subsequence' relation and hence we can apply Kruskal's Tree Theorem. It is not hard to see that if $\mathfrak T_{{\mathfrak {a}}_{M_1}}\leq\mathfrak T_{{\mathfrak {a}}_{M_2}}$, then ${\mathfrak {a}}_{M_1}\unlhd{\mathfrak {a}}_{M_2}$.
\endproof

\begin{defn}[eventuality; realization time]
Let
$\mathfrak D=\<W,R,g,t\>$
be a non-deterministic quasimodel.

An {\em eventuality} is any formula of the form $\neg\ast\psi\in {\mathrm{sub}}_\pm(\varphi)$.

Given a path $\vec w\in W^g$, $N\geq 0$ and an eventuality $\neg\ast\psi\in t(w_N)$, define the {\em realization time} of $\neg\ast\psi$ at $N$, denoted $\rho_N^{\neg\ast\psi}\(\vec w\)$, to be the least $K\geq N$ such that $\neg\ast\psi\in t(w_K)$. In case that no such $K$ exists set $\rho_N^{\neg\ast\psi}\(\vec w\)=\infty$.

Likewise, define
\[\rho_N\(\vec w\)=\cbra\rho_N^{\neg\ast\psi}\(\vec w\):\neg\ast\psi\in t(w_N)\cket.\]
Let $\rho^{\infty}_N\(\vec w\)$ be the maximum element of $\rho_N\(\vec w\)$ and $\rho^{<\infty}_N\(\vec w\)$ be the maximum finite element. In case that one of the sets being considered is empty, take zero instead of the maximum.
\end{defn}
\begin{defn}[efficiency]
Let
$\vec{\mathfrak a}=\<\mathfrak a_n\>$
be a finite or infinite path of local Kripke frames.

An {\em inefficiency} in $\vec{\mathfrak a}$ is a triple
$N\leq M_1<M_2$
such that ${\mathfrak {a}}_{M_1}\unlhd{\mathfrak {a}}_{M_2}$, $M_2<\rho^{\infty}_N(\vec{\mathfrak a})$ and
$\rho_N(\vec{\mathfrak a})\cap(M_1,M_2)=\varnothing.$

A finite or infinite path $\<\mathfrak a_n\>$ is {\em efficient} if it contains no inefficiencies and, for all $n\geq 0$,
$\|\mathfrak a_{n+1}\|\leq\|\mathfrak a_n\|+|\varphi|.$
\end{defn}

Roughly, the previous definition says that if the same state occurs twice in a row in an efficient path, some eventuality must have been realized in the middle. Otherwise, the path between them gives us a sort of loop which we could simply skip. Furthermore, efficiency gives us a way to guarantee that a path is realizing.

\begin{lem}\label{eff}
For all $\mathfrak a\in I_\omega(\varphi)$, there is a realizing path beginning on $\mathfrak a$ if and only if there is an efficient path beginning on $\mathfrak a$.
\end{lem}
\proof
First suppose that we have an efficient path
$\vec{\mathfrak a}=\<\mathfrak a_n\>_{n\geq 0}.$
We claim that $\vec{\mathfrak a}$ is realizing.

Let $N\geq 0$. We will show that $\rho^{\infty}_N\(\vec {\mathfrak a}\)<\infty$, and therefore that all eventualities of $\mathfrak a_N$ are realized.

By Lemma \ref{kruskapp}, there exist indices
$\rho^{<\infty}_N\(\vec{\mathfrak a}\)\leq M_1<M_2$
such that $\mathfrak a_{M_1}\unlhd \mathfrak a_{M_2}$. 

Since we know that $\vec{\mathfrak a}$ is efficient, this cannot produce an inefficiency. But no eventualities of $\mathfrak a_N$ occur between $M_1$ and $M_2$, so we must have $\rho^{\infty}_N\(\vec {\mathfrak a}\)\leq M_2$, and all eventualities of $\mathfrak a_N$ are thus realized by time $M_2$.

Since $N$ was arbitrary, it follows that the path is realizing.

For the other direction, suppose
$\vec{\mathfrak a}=\<\mathfrak a_n\>_{n\geq 0}$
is a realizing path. We claim that we can obtain an efficient path from $\vec{\mathfrak a}$ by removing all inefficiencies.

Let us first show how to remove a single inefficiency.

Suppose $\mathfrak a_{M_1}\unlhd\mathfrak a_{M_2}$ and this produces an inefficiency.

Then, clearly $\mathfrak a_{M_1}\rra\mathfrak a_{M_2+1}$, and we get a sequence
\[\mathfrak a_0\rra...\mathfrak a_{M_1}\rra\mathfrak a_{M_2+1}\rra...\]
which we can then reduce using Lemma \ref{reduction} to a sequence
\[\mathfrak a_0\rra\mathfrak a'_1\rra...\mathfrak a'_{n}\rra\mathfrak a'_{n+1}\rra...\]
with
$\|\mathfrak a'_{n+1}\|\leq\|\mathfrak a'_n\|+|\varphi|.$

Now, to obtain an efficient sequence, we can apply this process countably many times: first we ensure that no inefficient loops start at time $0$, then at time $1$, etc. The end result is well-defined because for all $n\geq 0$, the $n^{\mathrm{th}}$ element in the sequence stabilizes in finite time. This gives us an efficient realizing sequence.
\endproof

\begin{defn}[extension function]
Let
$\mathfrak D=\<W,R,g,t\>$
be a non-de\-ter\-mi\-nis\-tic quasimodel for a formula $\varphi$.

An {\em extension function} is a function
$\epsilon:W\to W^g,$ where $\epsilon(w)=\<\epsilon_n(w)\>_{n\geq 0}$ satisfies $\epsilon_0(w)=w$.
\end{defn}

Extension functions give us a canonical way to include points in realizing sequences. If we have an extension function on a typed Kripke frame, this gives us a way to guarantee that the transition relation is $\omega$-sensible.

\begin{defn}[family of paths]
A {\em family of paths} is a pair
$\mathfrak P=\<A,\epsilon\>,$
where $A\subseteq I_\omega(\varphi)$ is open and $\epsilon$ is an extension function assigning a realizing path in $A$ to each $\mathfrak a\in A$.

Likewise, a {\em partial family of paths of depth $N$} is a pair
$\mathfrak P^N=\<A^N,\epsilon^N\>,$
where $A\subseteq I_N(\varphi)$ and, for all $k\leq N$ and $\mathfrak a\in A\cap I_k(\varphi),$ $\epsilon^N(\mathfrak a)$ is a path of length $N-k+1$ in $A$ such that $\epsilon^N_0(\mathfrak a)=\mathfrak a$.

In either case,  $\mathfrak P$ is {\em efficient} if for all $\mathfrak a\in A$, $\epsilon(\mathfrak a)$ is efficient. $\mathfrak P$ {\em satisfies} $\varphi$ if there exists $\mathfrak a_\ast\in A\cap I_0(\varphi)$ such that $\varphi\in t(\mathfrak a_\ast)$.
\end{defn}
\begin{lem}\label{fam}
A formula $\varphi$ is satisfiable if and only if there exists an efficient family of paths $\mathfrak P$ satisfying $\varphi$.
\end{lem}

\proof
If $\mathfrak M$ is a model satisfying $\varphi$, then we can use Lemma \ref{eff} to assign an efficient path $\epsilon (\mathfrak a)$ to each $\mathfrak a\in \textrm{dom}(\chi^\ast)$.

This gives us an efficient family of paths
$\<\textrm{dom}(\chi^\ast),\epsilon\>.$

Conversely, it is easy to see that if we have an efficient family of paths, we also have a non-deterministic quasimodel; since $\epsilon$ gives us realizing paths in $A$, it follows that the relation $\rra$ is $\omega$-sensible on $\mathfrak I_\omega(\varphi)\upharpoonright A$.
\endproof

\begin{thm}
The set of all valid formulas of $\mathcal{DTL}$ is recursively enumerable.
\end{thm}
\proof
The strategy is to enumerate all efficient partial families of paths. This can be done, since there are only finitely many partial families of paths of any fixed depth $N$. We claim that $\varphi$ is valid if and only there exists $N\geq 0$ such that no efficient family of paths of depth $N$ satisfies $\neg\varphi$.

If $\neg\varphi$ is satisfiable, by Lemma \ref{fam}, there exists an efficient family of paths $\mathfrak P$ satisfying $\neg\varphi$. This immediately gives us an infinite sequence of partial families
$\mathfrak P^n=\mathfrak P\upharpoonright I_n(\varphi)$ satisfying $\neg\varphi$.

Conversely, if the search does not terminate, we can use K\"{o}nig's Lemma to find an increasing sequence
$\<\mathfrak P^n\>_{n\geq 0},$
satisfying $\neg\varphi$, where
$\mathfrak P^n=\mathfrak P^{n+1}\upharpoonright I_{n}(\varphi)$ for all $n$.

We can then define
$A=\bigcup_{n\geq 0}A^n,$ and for $K\geq 0$ and $\mathfrak a\in A^K$, set
$\epsilon_n(\mathfrak a)=\epsilon^{K+n}_n(\mathfrak a)$. Clearly the path $\epsilon(\mathfrak a)$ is efficient.

Thus,
$\mathfrak P=\<A,\epsilon\>$
gives us an efficient family of paths satisfying $\neg\varphi$, as desired.
\endproof

Unfortunately, the procedure we have just described does not suggest an obvious proof system for $\mathcal{DTL}$, and the above model-search algorithm is the only recursive enumeration of valid fomulas that we offer here. One axiomatization is suggested in $\cite{kremer}$; the question of its completeness remains open.

\bibliographystyle{plain}

\end{document}